\documentclass[12pt]{article}
\usepackage{a4wide,graphics,bm,color}
\usepackage{amsmath,latexsym,amssymb,amsfonts}
\usepackage{ijbc}
\usepackage{tabls}
\usepackage{array,delarray}
\usepackage{graphicx}
\begin{document}
\begin{titlepage}
  \title{CONTINUATION OF CONNECTING ORBITS IN 3D-ODES: (II)
    CYCLE-TO-CYCLE CONNECTIONS} \author{E.J. DOEDEL$^1$, B.W.
    KOOI$^2$, YU.A. KUZNETSOV$^3$, and G.A.K. van VOORN$^2$
    \\ \\
    {\it $^1$Department of Computer Science, Concordia University,}\\
    {\it 1455 Boulevard de Maisonneuve O., Montreal, Quebec, H3G 1M8, Canada}\\
    {\tt doedel@cs.concordia.ca}
    \\[1mm]
    {\it $^2$Department of Theoretical Biology,  Vrije Universiteit,}\\
    {\it de Boelelaan 1087, 1081 HV Amsterdam,
      the Netherlands}\\
    {\tt kooi@bio.vu.nl, george.van.voorn@falw.vu.nl}
    \\[1mm]
    {\it $^3$Department of Mathematics, Utrecht University}\\
    {\it Budapestlaan 6, 3584 CD Utrecht,
      the Netherlands}\\
    {\tt kuznet@math.uu.nl}
    \\[1mm]
    } \date{\today} \maketitle \abstract{In Part I of this paper we
    discussed new methods for the numerical continuation of
    point-to-cycle connecting orbits in 3-dimensional autonomous ODE's
    using projection boundary conditions. In this second part we extend
    the method to the numerical continuation of cycle-to-cycle
    connecting orbits. In our approach, the projection boundary
    conditions near the cycles are formulated using eigenfunctions
    of the associated adjoint variational equations, avoiding costly
    and numerically unstable computations of the monodromy matrices.
    The equations for the eigenfunctions are included in the defining
    boundary-value problem, allowing a straightforward implementation
    in {\sc auto}, in which only the standard features of the software
    are employed. Homotopy methods to find the connecting orbits are
    discussed in general and illustrated with an example from
    population dynamics. Complete {\sc auto} demos, which can be
    easily adapted to any autonomous 3-dimensional ODE system, 
    are freely available.}  
\end{titlepage}

\twocolumn
\section{Introduction}\label{sec:intro}

In a diversity of scientific fields bifurcation theory is used for the
analysis of systems of ordinary differential equations (ODE's) under
parameter variation. Many interesting phenomena in ODE systems are
linked to global bifurcations. Examples of such are overharvesting in
ecological models with bistability properties (Bazykin, 1998\nocite{Bazykin1998}; 
Antonovsky et al., 1990\nocite{Antonovskyetal1990}; Van Voorn et al.,
2007\nocite{Voornetal2007}), and the occurrence and disappearance of
chaotic behaviour in such models.  For example, it has been shown (see
Kuznetsov et al., 2001\nocite{Kuznetsovetal2001} and Boer et al.,
1999\nocite{Boeretal1999}, 2001\nocite{Boeretal2001}) that chaotic 
behaviour of the classical food chain models is associated with 
global bifurcations of point-to-point,
point-to-cycle, and cycle-to-cycle connecting orbits.

In Part I of this paper (Doedel et al., 2007\nocite{Doedeletal2007})
we discussed \textit{heteroclinic} connections between equilibria and
cycles. Here we look at connections that link a cycle to itself (a
\textit{homoclinic} cycle-to-cycle connection, for which the cycle is
necessarily saddle), or to another cycle (a \textit{heteroclinic}
cycle-to-cycle connection). Orbits homoclinic to the same hyperbolic
cycle are classical objects of the Dynamical Systems Theory. It is
known thanks to Poincare (1879)\nocite{Poicare1879}, Birkhoff
(1935)\nocite{Birkgoff1935}, Smale (1963)\nocite{Smale1963}, Neimark
(1967)\nocite{Neimark1967}, and L.P. Shilnikov
(1967)\nocite{Shilnikov1967} that a transversal intersection of the
stable and unstable invariant manifolds of the cycle along such an
orbit implies the existence of infinite number of saddle cycles
nearby. Disappearance of the intersection via collision of two
homoclinic orbits ({\em homoclinic tangency}) is an important global
bifurcation for which the famous H\'{e}non map turns to be a model
Poincar\'{e} mapping (Gavrilov and Shilnikov, 1972;
\nocite{GavrilovShilnikov72,GavrilovShilnikov73} Palis and Takens,
1993\nocite{PalisTakens1993}, see also Kuznetsov, 2004\nocite{Kuznetsov2004}).

Numerical methods for homoclinic orbits to equilibria have been
devised by Doedel and Kernevez (1986, but see Doedel et al.,
1997\nocite{Doedeletal1997}), who approximated homoclinic orbits by
periodic orbits of large but fixed period. Beyn
(1990\nocite{Beyn1990}) developed a direct numerical method for the
computation of such connecting orbits and their associated parameter
values, based on integral conditions and a truncated boundary value
problem (BVP) with projection boundary conditions.

The continuation of homoclinic connections in \textsc{auto} (Doedel et
al., 1997\nocite{Doedeletal1997}) improved with the development of
HomCont by Cham\-pneys and Kuz\-net\-sov
(1994\nocite{ChampneysKuznetsov1994}) and Cham\-pneys et
al. (1996\nocite{Champneysetal1996}). However, it is only suited for the
continuation of bifurcations of homoclinic point-to-point connections
and some heteroclinic point-to-point connections. A modification of
this software was introduced by Demmel et
al.~(2000\nocite{Demmeletal2000}), that uses the continuation of
invariant subspaces  (\textsc{CIS}-algorithm) for the location and 
continuation of homoclinic point-to-point connections.

Dieci and Rebaza (2004\nocite{DieciRebaza2004}) have also made
significant progress recently, by developing methods to continue
point-to-cycle and cycle-to-cycle connecting orbits based on another
work by Beyn (1994\nocite{Beyn1994}). Their me\-thod employs a
multiple shooting technique and requires the numerical solving for the
monodromy matrices associated with the periodic cycles involved in the
connection.

Our previous paper (Doedel et al., 2007\nocite{Doedeletal2007}) dealt
with a method for the detection and continuation of point-to-cycle
connections. Here this method is adapted for the continuation of
homoclinic and heteroclinic cycle-to-cycle connections. The me\-thod
is set up such that the homoclinic case is essentially a heteroclinic
case where the same periodic orbit (but not the periodic solution) is
doubled. In Section~\ref{sec:algorithm} we give a short overview of a
BVP formulation to solve a heteroclinic cycle-to-cycle problem. In
Section~\ref{sec:definingsystems} it is shown how boundary conditions
are implemented. In Section~\ref{sec:startingstrategies} we discuss
starting strategies to obtain approximate connecting orbits using
homotopy. In Section~\ref{sec:implementation} the BVP is made suitable
for numerical implementation.

Results are presented of the continuation of a homoclinic
cycle-to-cycle connection in the standard three-level food chain model
in Section~\ref{sec:results}. Boer et al.~(1999\nocite{Boeretal1999})
previously numerically obtained the two-parameter continuation curve
of this connecting orbit using a shooting method, combined with the
Poincar\'{e} map technique. In the previous part of this paper (Doedel
et al., 2007\nocite{Doedeletal2007}) we reproduced the results for the
structurally stable heteroclinic point-to-cycle connection of the same
food chain model using the homotopy method. In this paper we discuss
how the homoclinic cycle-to-cycle connection can be detected, and
continued in parameter space using the homotopy method. Also, it is
set up such that it can be used as well for a heteroclinic
cycle-to-cycle connection.

\section{Truncated BVPs with projection BCs}\label{sec:algorithm}

\begin{figure*}[htbp]
  \begin{center}
    \includegraphics[width=17.0cm]{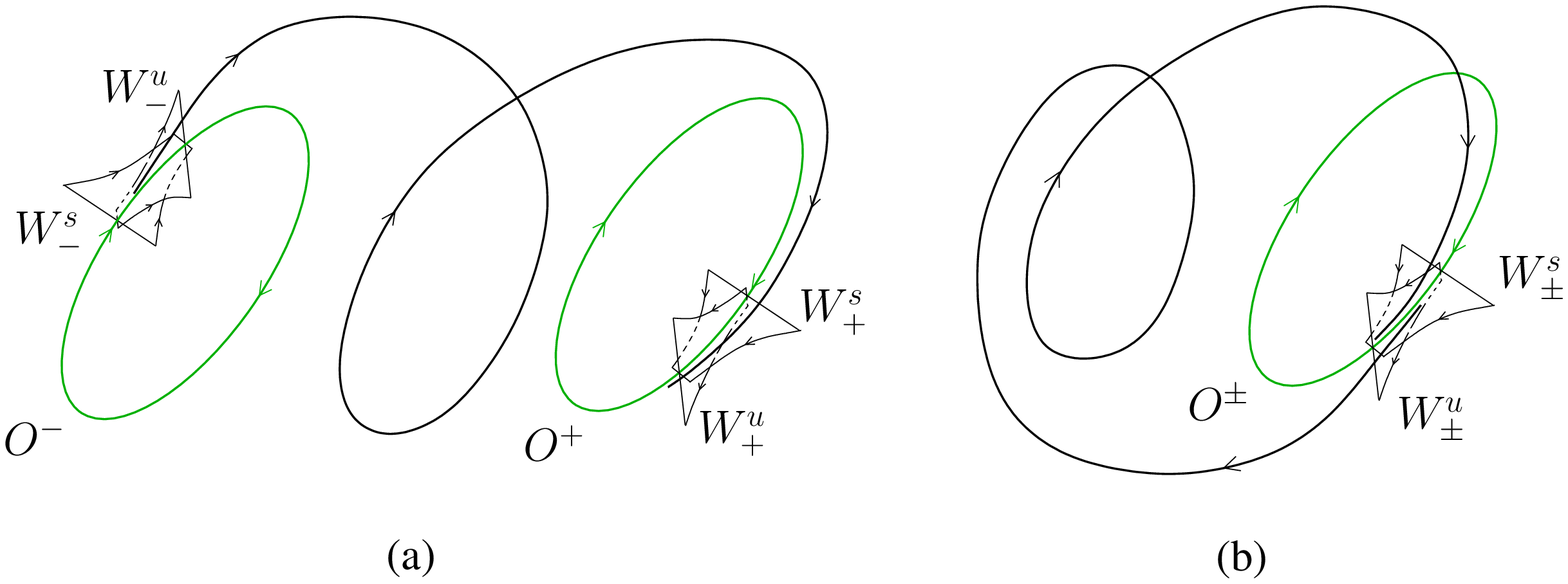}
  \end{center}
  \caption{Cycle-to-cycle connecting orbits in ${\mathbb R}^3$: (a) heterolinic
    orbit, $O^{+} \neq O^{-}$; (b) homoclinic orbit, $O^{+} = O^{-}$.}
  \label{fig:connections} 
\end{figure*}

Before presenting the BVP that describes a cycle-to-cycle connection,
let us first set up some notation. Consider a general system of ODEs
\begin{equation}
\label{eqn:ODE}
  \frac{du}{dt} = f(u,\alpha),
\end{equation}
where $f: \mathbb{R}^n \times \mathbb{R}^p \to \mathbb{R}^n$ is
sufficiently smooth, given that state variables $u\in \mathbb{R}^n$,
and control parameters $\alpha \in \mathbb{R}^p$. Thus, the dimension
of the state space is $n$ and the dimension of the parameter space is
$p$. The (local) flow generated by (\ref{eqn:ODE}) is denoted by
$\varphi^t$. Whenever possible, we will not indicate explicitly the 
dependence of various objects on parameters.

We assume that both $O^{-}$ and $O^{+}$ are saddle
limit cycles of (\ref{eqn:ODE}). A solution $u(t)$ of
(\ref{eqn:ODE}) for fixed $\alpha$ defines a {\em connecting orbit}
from $O^{-}$ to $O^{+}$ if
\begin{equation}
  \label{eqn:ASSYMPT}
  \lim_{t \to \pm \infty}{\rm dist}(u(t), O^{\pm}) = 0\;.
\end{equation}
(Figure~\ref{fig:connections} depicts such a connecting orbit in the
3D-space.)  Since $u(t + \tau)$ satisfies (\ref{eqn:ODE}) and
(\ref{eqn:ASSYMPT}) for any phase shift $\tau$, an additional phase
condition
\begin{equation}
  \label{eqn:PHASE}
  \psi[u,\alpha] = 0 \;,
\end{equation}
should be imposed to ensure uniqueness of the connecting
solution. This condition will be specified later.

\begin{table}[htbp]
  \caption[]{\label{tab:values}{\protect\small List of notation used
  in the paper.}}
  \begin{center}
    \begin{tabular}{p{0.7cm}p{6.5cm}}\hline      
      sym. & meaning\\
      \hline
      $x^{\pm}$ & Periodic solution\\
      $v^{\pm}$ & Eigenfunction\\
      $w^{\pm}$ & Scaled adjoint eigenfunction\\
      $u$ & Connection\\
      $\alpha$ & Bifurcation parameters\\
      $O^{+}$ & Limit cycle where connection ends\\
      $O^{-}$ & Limit cycle where connection starts\\
      $W^s_{+}$ & Stable manifold of the cycle $O^{+}$\\
      $W^u_{-}$ & Unstable manifold of the cycle $O^{-}$\\
      $\mu^+_u$ & Unstable multiplier of the cycle $O^{+}$\\
      $\mu^-_s$ & Stable multiplier of the cycle $O^{-}$\\
      $\mu^-_u$ & Unstable multiplier of the cycle $O^{-}$\\
      $\mu^{+}$ & Adjoint multiplier $1/\mu^{+}_u$\\
      $\mu^{-}$ & Adjoint multiplier $1/\mu^{-}_s$\\
      $\lambda^{\pm}$ & $\ln(\mu^{\pm})$\\
      $T^{\pm}$ & Period of the cycle $O^{\pm}$\\
      $T$ & Connection time\\
      \hline
    \end{tabular}
  \end{center}
\end{table}

For numerical approximations, the asymptotic conditions
(\ref{eqn:ASSYMPT}) are substituted by {\em projection boundary
conditions} at the end-points of a large {\em truncation interval}
$[\tau_{-},\tau_{+}]$, following Beyn (1994\nocite{Beyn1994}). It is
prescribed that the points $u(\tau_{-})$ and
$u(\tau_{+})$ belong to the linear subspaces, which are tangent
to the unstable and stable invariant manifolds of $O^{-}$ and
$O^{+}$, respectively.

Now, denote by $x^{\pm}(t)$ a periodic solution (with minimal
period $T^{\pm}$) corresponding to $O^{\pm}$ and
introduce the {\em monodromy matrix}
$$
M^{\pm}=\left.D_x \varphi^{T^{\pm}}(x)\right|_{x=x^{\pm}(0)},
$$ 
{\it i.e.} the linearization matrix of the $T^{\pm}$-shift along
orbits of (\ref{eqn:ODE}) at point $x_0^{\pm} =
x^{\pm}(0) \in O^{\pm}$. Its eigenvalues 
are called {\em Floquet multipliers}, of which one (trivial) equals
1. Let $m_s^{+} = n_s^{+} + 1$ be the dimension of the stable
invariant manifold $W^s_{+}$ of the cycle
$O^{+}$, where $n_s^{+}$ is the number of its multipliers
satisfying
$$
|\mu| < 1.
$$ 
Along the same line, $m_u^{-} = n_u^{-} + 1$ is the dimension of the
unstable invariant manifold $W^u_{-}$ of the cycle
$O^{-}$, and $n_u^{-}$ is the number of its multipliers
satisfying
$$
|\mu| > 1.
$$

To have an isolated branch of cycle-to-cycle connecting orbits of
(\ref{eqn:ODE}) it is necessary that
\begin{equation}
  \label{eqn:FREEPARS1}
  p = n-m_s^{+} - m_u^{-} + 2 \;,
\end{equation}
(see Beyn, 1994\nocite{Beyn1994}).

The projection boundary conditions in this case become
\begin{equation}\label{eqn:BVPc2c}
    L^{\pm} (u(\tau_{\pm}) - x^{\pm}(0))  = 0 \;,
\end{equation}
where $L^{-}$ is a $(n-m^{-}_u) \times n$ matrix whose rows
form a basis in the orthogonal complement to the linear subspace that
is tangent to $W^u_{-}$ at $x^{-}(0)$. Similarly,
$L^{+}$ is a $(n-m^{+}_s) \times n$ matrix, such that its rows
form a basis in the orthogonal complement to the linear subspace that
is tangent to $W^s_{+}$ at $x^{+}(0)$.

The above construction also applies in the case when $O^{+}$
and $O^{-}$ coincide, {\it i.e.} we deal with a {\em homoclinic
orbit} to a saddle limit cycle $O^{+} =
O^{-}$. Note that, in general, the base points
$x^{\pm}(0) \in O^{\pm}$ remain different (and so do
the periodic solutions $x^{\pm}(t)$).

It can be proved that, generically, the truncated BVP composed of
(\ref{eqn:ODE}), a truncation of (\ref{eqn:PHASE}), and
(\ref{eqn:BVPc2c}), has a unique solution branch $(\hat{u} (t,
\hat{\alpha}), \hat{\alpha})$, provided that (\ref{eqn:ODE}) has a
connecting solution branch satisfying (\ref{eqn:PHASE}) and
(\ref{eqn:FREEPARS1}).

The truncation to the finite interval $[\tau_{-},\tau_{+}]$ causes an
error. If $u$ is a generic connecting solution to
(\ref{eqn:ODE}) at parameter $\alpha$, then the following estimate
holds in both cases:
\begin{equation}
  \label{ERROR}
  \|(u|_{[\tau_{-},\tau_{+}]},\alpha) - (\hat{u},\hat{\alpha}) \|\leq
    C {\rm e}^ {-2\min(\mu_{-}|\tau_{-}|,\mu_{+}|\tau_{+}|)},
\end{equation}
where $\|\cdot\|$ is an appropriate norm in the space
$C^1([\tau_{-},\tau_{+}], {\mathbb R}^n) \times {\mathbb R}^p$,
$u|_{[\tau_{-},\tau_{+}]}$ is the restriction of $u$ to the truncation
interval, and $\mu_{\pm}$ are determined by the eigenvalues of the
monodromy matrices. For exact formulations, proofs, and references to
earlier contributions, see Pampel (2001\nocite{Pampel2001}) and Dieci
and Rebaza (2004\nocite{DieciRebaza2004a,DieciRebaza2004}, including Erratum).

\section{New defining systems in $\mathbb{R}^3$}\label{sec:definingsystems}

In this section we show how to implement the boundary conditions
(\ref{eqn:BVPc2c}). We consider the case $n=3$ where $O^-$ and $O^+$
are saddle cycles and therefore always $m_s^{-}=m_u^{-}=2$ and
$m_s^{+}=m_u^{+}=2$. Substitution in (\ref{eqn:FREEPARS1}) gives the
number of free parameters for the continuation $p=1$. 


Note that the complete BVP will consist of 15 equations (2 saddle
cycles, 2 eigendata for these cycles, and the connecting orbit) and 19
boundary conditions.

\begin{figure*}[ht]
  \begin{center}
    \includegraphics[width=14.0cm]{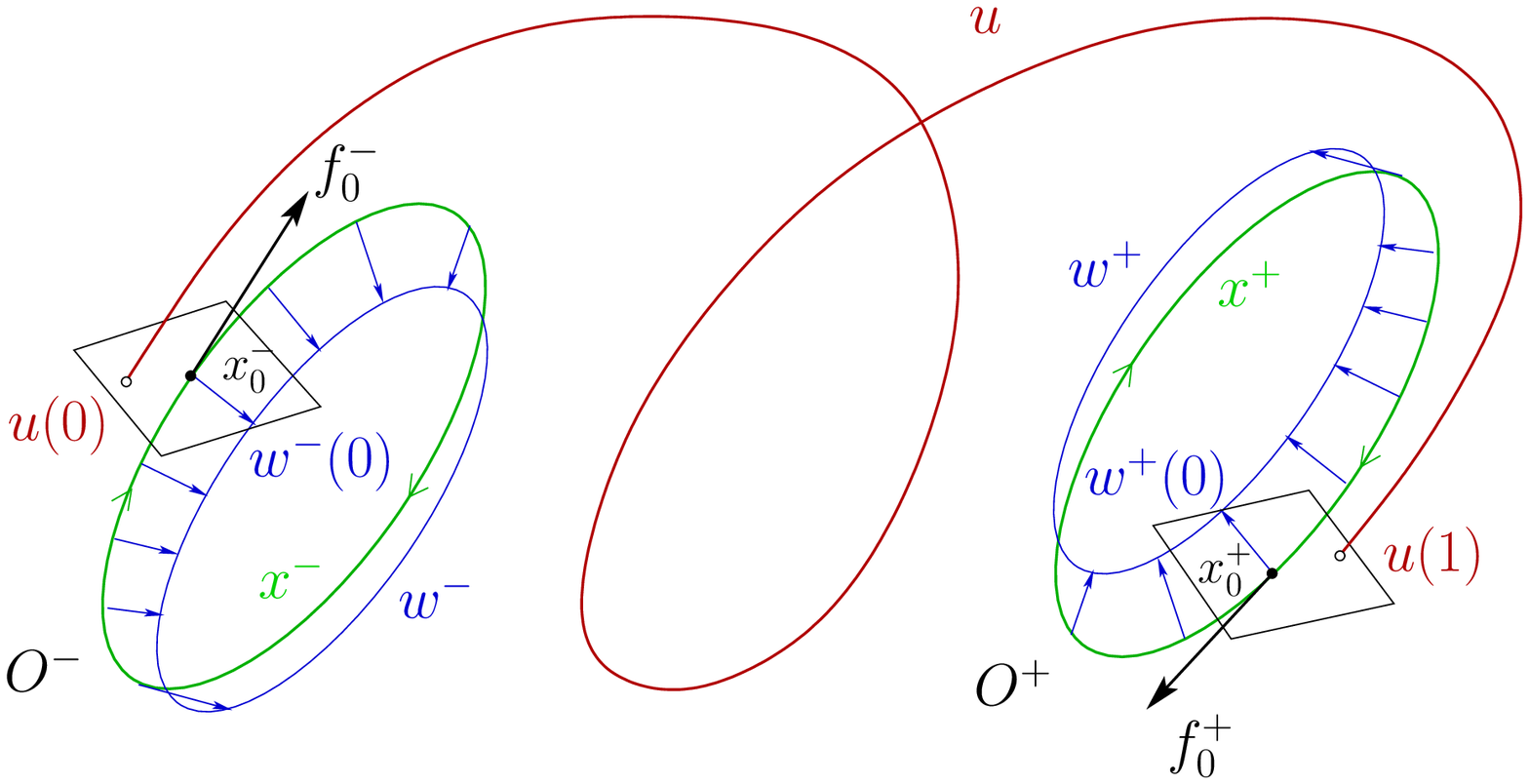}
  \end{center}
  \caption{Ingredients of a BVP to approximate a heteroclinic connecting orbit. The
  homoclinic cycle-to-cycle connection is also approached as the
  heteroclinic case, where two saddle cycles coincide.}
  \label{fig:conBVP} 
\end{figure*}

\subsection{The cycle and eigenfunctions}

To compute the saddle limit cycles $O^-$ and $O^+$ involved in the
heteroclinic connection (see Figure~\ref{fig:conBVP}) we need a BVP.
The standard periodic BVP can be used
\begin{equation}
\label{eqn:L}
    \left\{
      \begin{array}{rcl}
        \dot{x}^{\pm} - f(x^{\pm},\alpha) & = & 0 \;,\\
         x^{\pm}(0) - x^{\pm}(T^{\pm}) & =  & 0\;, 
      \end{array}
    \right.  
\end{equation}
A unique solution of this BVP is determined by using an appropriate
phase condition, which is actually a boundary condition for the
truncated connecting solution, and which will be introduced below.

To set up the projection boundary condition for the truncated
connecting solution $u$ near $O^{\pm}$, we also need a vector, say
$w^{+}(0)$, that is orthogonal at $x^{+}(0)$ to the stable
manifold $W^s_{+}$ of the saddle limit cycle $O^{+}$, as well another 
vector, say $w^{-}(0)$, that is orthogonal at $x^{-}(0)$ to the unstable
manifold $W^u_{-}$ of the saddle limit cycle $O^{-}$ (see
Figure~\ref{fig:conBVP}). Each vector $w^{\pm}(0)$ can be obtained from
an {\em eigenfunction} $w^{\pm}(t)$ of the {\em adjoint variational
problem} associated with (\ref{eqn:L}), corresponding to 
eigenvalue $\mu^{\pm}$. These eigenvalues satisfy
$$
\mu^{+} = \frac{1}{\mu^{+}_u} \;,\; \mu^{-} = \frac{1}{\mu^{-}_s}\;,
$$ 
where $\mu^{+}_u$ and $\mu^{-}_s$ are the multipliers of the monodromy matrix
$M^{\pm}$ with
$$
|\mu^{+}_u| > 1 \; , \; |\mu^{-}_s| < 1\;.
$$
The corresponding BVP is
\begin{equation}
\label{eqn:EF_original}
\left\{
      \begin{array}{rcl}
        \dot{w}^{\pm} + f_u^{\rm T}(x^{\pm},\alpha) w^{\pm} & = & 0 \;,\\
        w^{\pm}(T^{\pm}) - \mu^{\pm} w^{\pm}(0) & = & 0\;,\\
        \langle w^{\pm}(0),w^{\pm}(0)  \rangle - 1& = &0 \;,
      \end{array}
    \right.
\end{equation}
where $x^{\pm}$ is the solution of (\ref{eqn:L}).  In our
implementation the above BVP is replaced by an equivalent BVP
\begin{equation}
\label{eqn:EF}
\left\{
      \begin{array}{rcl}
        \dot{w}^{\pm} + f_u^{\rm T}(x^{\pm},\alpha) w^{\pm} + \lambda^{\pm}w^{\pm}& = & 0 \;,\\
        w^{\pm}(T^{\pm}) - s^{\pm} w^{\pm}(0) & = & 0\;,\\ 
        \langle w^{\pm}(0), w^{\pm}(0) \rangle - 1& = & 0 \;,
      \end{array}
    \right.
\end{equation}
where $s^{\pm} = {\rm sign\;} \mu^{\pm}$ and 
$$
\lambda^{\pm} = \ln|\mu^{\pm}|\;.
$$
(See Appendix of Part I, Doedel et al., 2007\nocite{Doedeletal2007}).

In (\ref{eqn:EF}), the boundary conditions become periodic or
anti-periodic, depending on the sign of the multiplier $\mu^{\pm}$,
while the logarithm of its absolute value appears in the variational
equation. This ensures high numerical robustness.

Given $w^{\pm}$ satisfies (\ref{eqn:EF}), the projection boundary
conditions (\ref{eqn:BVPc2c}) become
\begin{equation}
  \label{eqn:PROJCYCLE}
  \langle w^{\pm}(0), u(\tau_{\pm}) - x^{\pm}(0) \rangle = 0.
\end{equation}

\subsection{The connection}

We use the following BVP for the connecting solution:
\begin{equation}
  \label{eqn:C}
  \left\{
  \begin{array}{rcl}
    \dot{u} - f(u,\alpha) & = & 0 \;,\\
    \langle f(x^{\pm}(0),\alpha), u(\tau_{\pm}) - x^{\pm}(0) \rangle & = & 0 \;.
  \end{array}
  \right.
\end{equation}
For each cycle, a phase condition is needed to select a unique periodic 
solution among those which satisfy (\ref{eqn:L}), {\it i.e.} to fix a {\em base
point} $x_0^{\pm}=x^{\pm}(0)$ on the cycle $O^{\pm}$ (see
Figure~\ref{fig:conBVP}).  For this we
require the end-point of the connection to belong to a plane
orthogonal to the vector $f(x^+(0),\alpha)$, and the starting point of
the connection to belong to a plane orthogonal to the vector
$f(x^-(0),\alpha)$. This allows the base points $x^{\pm}(0)$ to 
move freely and independently upon each other 
along the corresponding cycles $O^{\pm}$.

\subsection{The complete BVP}\label{Sec:CompleteBVP}

The complete truncated BVP to be solved numerically consists of
\begin{subequations}
  \label{eqn:BVP}
  \begin{align}
    \dot{x}^{\pm} - T^{\pm} f(x^{\pm},\alpha) & = 0, \label{BVP-2}\\
    x^{\pm}(0) - x^{\pm}(1) & = 0 \;, \label{BVP-3}\\
    \dot{w}^{\pm} + T^{\pm} f_u^{\rm T}(x^{\pm},\alpha) w^{\pm} + \lambda^{\pm}
    w^{\pm} & = 0 \;, \label{BVP-5}\\
    w^{\pm}(1) - s^{\pm} w^{\pm}(0) & =  0, \label{BVP-6}\\
    \langle w^{\pm}(0), w^{\pm}(0)  \rangle - 1& =  0\;, \label{BVP-7}\\
    \dot{u} - T f(u,\alpha) & = 0\;, \label{BVP-8}\\
    \langle f(x^{+}(0),\alpha), u(1) - x^{+}(0) \rangle & = 0\;,
    \label{BVP-9a}\\
    \langle f(x^{-}(0),\alpha), u(0) - x^{-}(0) \rangle & = 0\;,
    \label{BVP-9b}\\
    \langle w^{+}(0), u(1) - x^{+}(0) \rangle &= 0 \;, \label{BVP-4p}\\
    \langle w^{-}(0), u(0) - x^{-}(0) \rangle &= 0 \;, \label{BVP-4q}\\
    \| u(0) - x^{-}(0) \|^2 - \varepsilon^2 & = 0\;,
    \label{BVP-0} 
  \end{align}
\end{subequations}
where the last equation places the starting point $u(0)$ of the connection 
at a small fixed distance $\varepsilon > 0$ from the base point $x^{-}(0)$.  
The time variable is scaled to the unit interval $[0,1]$, so that both the
cycle periods $T^{\pm}$ and the connection time $T$ become parameters.
Hence, besides a component of $\alpha$, there are five more parameters
available for continuation: the connection time $T$, the cycle periods
$T^{\pm}$, and the multipliers $\lambda^{\pm}$.

\section{Starting strategies}\label{sec:startingstrategies}

The BVP described in the previous section are only usable when good
initial starting data are available. Usually, such data are not
present. Here we demonstrate how initial data can be generated through
a series of successive continuations in {\sc auto}, a method referred
to as {\em homotopy method}, first introduced by Doe\-del, Fried\-man
and Mon\-tei\-ro (1994)\nocite{Doedeletal1994} for point-to-point
problems and extended to point-to-cycle problems in Part I of this paper.

\subsection{Saddle cycles}

The easiest way to obtain the limit saddle cycles $O^{\pm}$, first
calculate a stable equilibrium using software like {\sc maple}, {\sc
matlab} or {\sc mathematica}. Then, using {\sc auto}, continue this
equilibrium up to an Andronov-Hopf bifurcation, where a stable limit cycle is
generated. A continuation of this cycle can
result in the detection of a fold bifurcation for the limit
cycle. This will yield a saddle limit cycle.

\subsection{Eigenfunctions}\label{sec:eigfunc}

In order to obtain an initial starting point for the connecting orbit
we require knowledge about the unstable manifold of the saddle limit
cycle $O^{-}$. Also, we need the linearized adjoint ``manifolds'' to
understand how the connecting orbit leaves $O^{-}$ and approaches $O^{+}$
(or the same cycle in the homoclinic case). For this,
we look at the {\em eigendata}. 

First consider the periodic BVP for $O^{-}$, 
\begin{equation}
\label{eqn:L-}
    \left\{
      \begin{array}{rcl}
        \dot{x}^{-} - T^{-}f(x^{-},\alpha) & = & 0 \;,\\
         x^{-}(0) - x^{-}(1) & =  & 0\;, 
      \end{array}
    \right.  
\end{equation}
to which we add the standard integral phase condition
\begin{equation}
\label{eqn:CYCLEPHASE}
\int_0^1 \langle \dot{x}^{-}_{old}(\tau), x^{-}(\tau) \rangle = 0\;,
\end{equation}
as well as a BVP similar to (\ref{eqn:EF_original}), namely
\begin{equation}
\label{eqn:EF_original_h}
\left\{
      \begin{array}{rcl}
        \dot{v} - T^{-} f_u(x^{-},\alpha) v & = & 0 \;,\\
        v(1) - \mu v(0) & = & 0\;,\\
        \langle v(0),v(0)  \rangle - h& = &0 \;.
      \end{array}
    \right.
\end{equation}
In (\ref{eqn:CYCLEPHASE}), $x^{-}_{old}$ is a reference periodic solution,
{\it e.g.} from the preceding continuation step.  The parameter
$h$ in (\ref{eqn:EF_original_h}) is a {\em homotopy parameter}, that
is set to zero initially. Then, (\ref{eqn:EF_original_h}) has a
trivial solution
$$
v(t) \equiv 0,\ \ h=0,
$$
for any real $\mu$. This family of the trivial solutions parametrized
by $\mu$ can be continued in {\sc auto} using a BVP consisting of
(\ref{eqn:L-}), (\ref{eqn:CYCLEPHASE}),
and (\ref{eqn:EF_original_h}) with free parameters $(\mu,h)$ and fixed
$\alpha$. The unstable Floquet multiplier of $O^{-}$ then corresponds
to a branch point at $\mu=\mu^{-}_u$ along this trivial solution family. 
{\sc auto} can accurately locate such a point and switch to
the nontrivial branch that emanates from it. This secondary family is
continued in $(\mu,h)$ until the value $h=1$ is reached, which gives a
normalized {\em eigenfunction} $v^{-}$ corresponding to the multiplier $\mu^{-}_u$.
Note that in this continuation the value of $\mu$ remains constant,
$\mu \equiv \mu^{-}_u$, up to numerical accuracy. For the initial starting
point of the connection we use $v^{-}(0)$.

The same method is applicable to obtain the nontrivial {\em scaled 
adjoint eigenfunctions} $w^{\pm}$ of the saddle cycles. For this, the BVP
\begin{equation}
  \label{eqn:EF_h}
  \left\{
  \begin{array}{rcl}
    \dot{w}^{\pm} + T^{\pm} f_u^{\rm T}(x^{\pm},\alpha) w^{\pm} + \lambda^{\pm}
    w^{\pm} & = & 0 \;,\\
    w^{\pm}(1) - s^{\pm} w^{\pm}(0) & = & 0\;,\\
    \langle w^{\pm}(0),w^{\pm}(0)  \rangle - h^{\pm}& = &0 \;,
  \end{array}
  \right.
\end{equation}
where $s^{\pm} = {\rm sign}(\mu^{\pm})$, replaces
(\ref{eqn:EF_original_h}). A branch point at $\lambda_1^{\pm}$ then
corresponds to the adjoint multiplier $\sigma^{\pm} {\rm
e}^{\lambda_1^{\pm}}$. After branch switching the desired eigendata
can be obtained.

\subsection{The connection}

Time-integration of (\ref{eqn:ODE}), in {\sc matlab} for instance,
can yield an initial connecting orbit, however, this only applies for
non-stiff systems. Nevertheless, mostly when starting sufficiently
close to the exact connecting orbit {\em in parameter space} the
method of {\em successive continuation} (Doe\-del, Fried\-man and
Mon\-tei\-ro, 1994\nocite{Doedeletal1994}) can be used to obtain an
initial connecting orbit.

Let us introduce a BVP that is a modified version of (\ref{eqn:BVP})
\begin{subequations}\label{eqn:BVP-homotopy}
  \begin{align}
    \dot{x}^{\pm} - T^{\pm} f(x^{\pm},\alpha) & = 0 \;, \label{BVP22a1}\\
    x^{\pm}(0) - x^{\pm}(1) & = 0 \;, \label{BVP22b1}\\
    \Phi^{\pm}[x^{\pm}] &= 0 \;, \label{BVP22c1}\\
    \dot{w}^{\pm} + T^{\pm} f_u^{\rm T}(x^{\pm},\alpha) w^{\pm} +
    \lambda^{\pm} w^{\pm} & = 0 \;, 
    \label{BVP22d}\\
    w^{\pm}(1) - s^{\pm} w^{\pm}(0) & =  0 \;, \label{BVP22e1}\\
    \langle w^{\pm}(0), w^{\pm}(0)  \rangle - 1& =  0 \;, \label{BVP22f1}\\
    \dot{u} - T f(u,\alpha) & = 0, \label{BVP22g}\\
    \langle f(x^{+}(0),\alpha), u(1) - x^{+}(0) \rangle 
    - h_{11} & =  0 \;, \label{BVP22h1}\\
    \langle f(x^{-}(0),\alpha), u(0) - x^{-}(0) \rangle 
    - h_{12} & =  0 \;, \label{BVP22h2}\\
    \langle w^{+}(0), u(1) - x^{+}(0) \rangle - h_{21} & = 0 \;, 
    \label{BVP24c1}\\
    \langle w^{-}(0), u(0) - x^{-}(0) \rangle - h_{22} & = 0 \;, 
    \label{BVP24c2}
  \end{align}
\end{subequations}
where each $\Phi^{\pm}$ in (\ref{BVP22c1}) defines any phase condition fixing the
base point $x^{\pm}(0)$ on the cycle $O^{\pm}$. An example of such a phase
condition is
$$
\Phi^{+}[x] = x_j(0) - a_j \;,
$$
where $a_j$ is the $j$th-coordinate of the base point of $O^{+}$ at some given
parameter values. Furthermore, $h_{jk}$, $j,k = 1,2,$ in (\ref{BVP22h1})--(\ref{BVP24c2}) 
are {\em homotopy parameters}.

For the approximate connecting orbit a
small step $\varepsilon$ is made in the direction of the unstable
eigenfunction $v^{-}$ of the cycle $O^{-}$: 
\begin{equation}
  \label{eqn:INIT1}
  u(\tau) = x^{-}(0) + \varepsilon v^{-}(0) {\rm e}^{\mu^{-}_u T^{-}\tau},~~
  \tau \in [0,1] \;,
\end{equation}
which provides an approximation to a solution of $\dot{u}=T^{-}f(u,\alpha)$ 
in the unstable manifold $W^u_{-}$ near $O^{-}$. 
After collection of the cycle-related data, eigendata
and the time-integrated approximated orbit, a solution to the above BVP can be 
continued in $(T,h_{11})$ and $(T,h_{12})$ for fixed value of $\alpha$
in order to make $h_{11} = h_{12} = 0$, while $u(1)$ is
near the cycle $O^+$, so that $T$ becomes sufficiently large. In the
next step, we then try to make $h_{21} = h_{22} = 0$, after which a
good approximate initial connecting orbit is obtained.

This solution is now used to activate one of the system parameters, say $\alpha_1$, 
and to continue a solution to the primary BVP (\ref{eqn:BVP}). 
Then, if necessary after having improved the connection first by a
continuation in $T$, continuation in $(\alpha_1,T)$ can be done to
detect limit points, using the standard fold-detection facilities of
{\sc auto}. Subsequently a fold curve can be continued in two
parameters, say $(\alpha_1,\alpha_2)$, for fixed $T$ using the
standard fold-following facilities in {\sc auto}.

\section{Implementation in AUTO}\label{sec:implementation}

Our algorithms have been implemented in \textsc{auto}, which solves
the boundary value problems using superconvergent {\em orthogonal
collocation} with adaptive meshes.  {\sc auto} can compute paths of
solutions to boundary value problems with integral constraints and
non-separated boundary conditions:
\begin{subequations}
  \label{eqn:AUTO}
  \begin{align}
    \dot{U}(\tau) - F(U(\tau),\beta) &=0\; , \; \; \tau \in [0,1],
    \label{5.1} \\
    b(U(0),U(1),\beta) &= 0\;,  \label{5.2} \\
    \int^1_0 q(U(\tau),\beta) d\tau &= 0\;, \label{5.3} 
  \end{align}
\end{subequations}
where
$$
U(\cdot),F(\cdot,\cdot) \in  {\mathbb R}^{n_d}, ~
b(\cdot,\cdot) \in {\mathbb R}^{n_{bc}}, ~
q(\cdot,\cdot) \in {\mathbb R}^{n_{ic}}, 
$$
and
$$
\beta \in {\mathbb R}^{n_{fp}}, 
$$
as $n_{fp}$ {\it free} parameters $\beta$ are allowed to vary, where
\begin{equation}
n_{fp} = n_{bc} + n_{ic}  - n_d + 1\;.
\label{5.4}
\end{equation}
The function $q$ can also depend on $F$, the derivative of $U$ with
respect to pseudo-arclength, and on $\hat{U}$, the value of $U$ at the
previously computed point on the solution family.

For our primary BVP problem (\ref{eqn:BVP}) in three dimensions we have
$$
n_d = 15, ~~ n_{ic}=0,
$$
and $n_{bc}=19$, so that any 5 free parameters are allowed to vary.

\section{Example: food chain model}\label{sec:results}

In this section we describe the performance of the BVP-method for the
detection and continuation of a cycle-to-cycle connecting orbit in the
standard food chain model, also used in Part I of this paper.

\subsection{The model}
The three-level food chain model from theoretical biology, based on
the Rosen\-zweig-MacAr\-thur (1963\nocite{RosenzweigMacArthur1963})
prey-predator model, is given by the following equations
\begin{equation}\label{eqn:RM3D}
  \left\{
  \begin{array}{rcl}
    \dot{x_1} & = & x_1 (1 - x_1) - f_1(x_1,x_2)\;,\\
    \dot{x_2} & = & f_1(x_1,x_2) - d_1 x_2 - f_2(x_2,x_3)\;,\\
    \dot{x_3} & = & f_2(x_2,x_3) - d_2 x_3\;,
  \end{array}
  \right.
\end{equation}
with Holling Type-II functional responses
$$
f_1(x_1,x_2) = \frac{a_1 x_1 x_2}{1 + b_1 x_1}\;
$$
and
$$
f_2(x_2,x_3) = \frac{a_2 x_2 x_3}{1 + b_2 x_2}\;.
$$ This standard model has been studied by several authors, see {\it
e.g.} Kuz\-net\-sov and Ri\-nal\-di
(1996\nocite{KuznetsovRinaldi1996}) and Kuz\-net\-sov et
al. (2001\nocite{Kuznetsovetal2001}) and references there.

\begin{figure}[htbp]
  \scalebox{0.7}{
    \includegraphics{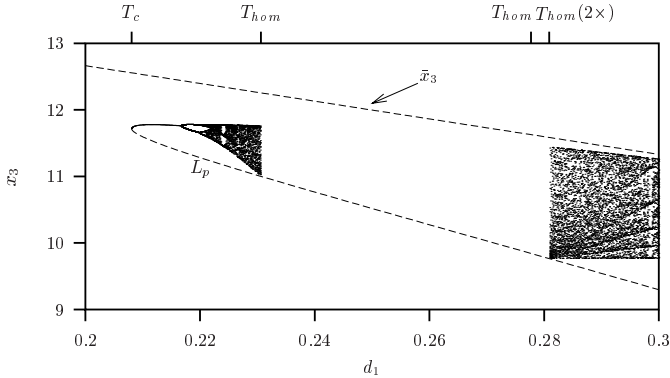}}
  \vspace{-1.0cm}
  \caption{
    One-parameter bifurcation diagram for $d_2 = 0.0125$. The
    equilibrium is indicated as $\bar{x}_3$. The dashed line $L_p$ is
    the $x_3$-value of the local minimum of an unstable (saddle) limit
    cycle.  At the point $T_c$ this limit cycle coincides with a
    stable limit cycle, of which the local minimum of $x_3$ is also
    shown. The stable limit cycle undergoes period doublings until
    chaos (the dense regions) is reached.  The two dense regions are
    separated by a region where homoclinic cycle-to-cycle connections
    to the limit cycle $L_p$ exist.  Both chaotic regions are bounded
    by a limit point of the homoclinic connection, indicated by
    $T_{hom}$.  Observe that near the right chaotic region, between
    two limit points exist secondary connecting orbits to the cycle.
    One of these limit points coincides with the limit point of the
    primary connecting orbit that forms the boundary of the region of
    chaos (hence $T_{hom}$ twice); after Boer et al., 1999).}
  \label{fig:fig8ch4}
\end{figure}
\begin{figure}
  \scalebox{0.7}{
    \includegraphics{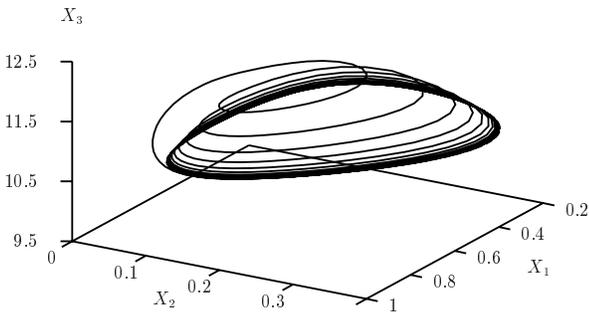}}
  \vspace{-1.0cm}
  \caption{Phase plot of the homoclinic cycle-to-cycle connection.}
  \label{fig:phplotrm}
\end{figure}

The death rates $d_1$ and $d_2$ are often used as bifurcation
parameters $\alpha_1$ and $\alpha_2$, respectively, with the other parameters set at $a_1 = 5$, $a_2 = 0.1$,
$b_1 = 3$, and $b_2 = 2$. For these parameter values the model
displays chaotic behaviour in a given parameter range of $d_1$ and
$d_2$ (Hastings and Powell, 1991\nocite{HastingsPowell1991}; Klebanoff
and Hastings, 1994\nocite{KlebanoffHastings1994}; McCann and Yodzis,
1995\nocite{McCannYodzis1995}). The region of chaos can be found starting from a
fold bifurcation at for instance $d_1 \approx 0.2080452$, $d_2 =
0.0125$, where two limit cycles appear. The stable branch then
undergoes a cascade of period-doublings (see Figure~\ref{fig:fig8ch4})
until a region of chaos is reached.

Previous work by Boer et al.~(1999\nocite{Boeretal1999},
2001\nocite{Boeretal2001}) has shown that the parameter region where
chaos occurs is intersected by homoclinic and heteroclinic global
connections, and that this region is partly bounded by a homoclinic
cycle-to-cycle connection, as shown in Figure~\ref{fig:fig8ch4}. These
results were obtained numerically using multiple shooting.

\subsection{Homotopy}

Using the technique discussed in this paper we first find the saddle
limit cycle for $d_1 = 0.25$, $d_2 = 0.0125$. Since the cycle $O$ is
both $O^{+}$ and $O^{-}$, we use the same initial base point
$$
x^{\pm}(0) = (0.839783,0.125284,10.55288)
$$
and the period $T^{\pm} = 24.28225$. The logarithms of the nontrivial adjoint
multipliers are
$$
\lambda^{+} = -0.4399607\;,\; \lambda^{-} = 6.414681\;.
$$
The starting point of the initial ``connecting'' orbit is calculated
by taking the base point $x^{-}(0)$ and multiplying the eigenfunction $v^{-}(0)$
by $\varepsilon = -0.001$
\begin{equation}
  u(0) = x^{-}(0) + \varepsilon v^{-}(0) \;,
\end{equation}
where
$$
v^{-}(0) = (-1.5855 \cdot 10^{-2}, 2.6935 \cdot 10^{-2},-0.99951)\;,
$$
and the resulting
$$
u(0) = (0.839789,0.125274,10.55324)\;.
$$
The connection time $T = 503.168$.

\begin{figure*}\hspace{0.cm}
  \scalebox{0.7}{
  { { \includegraphics{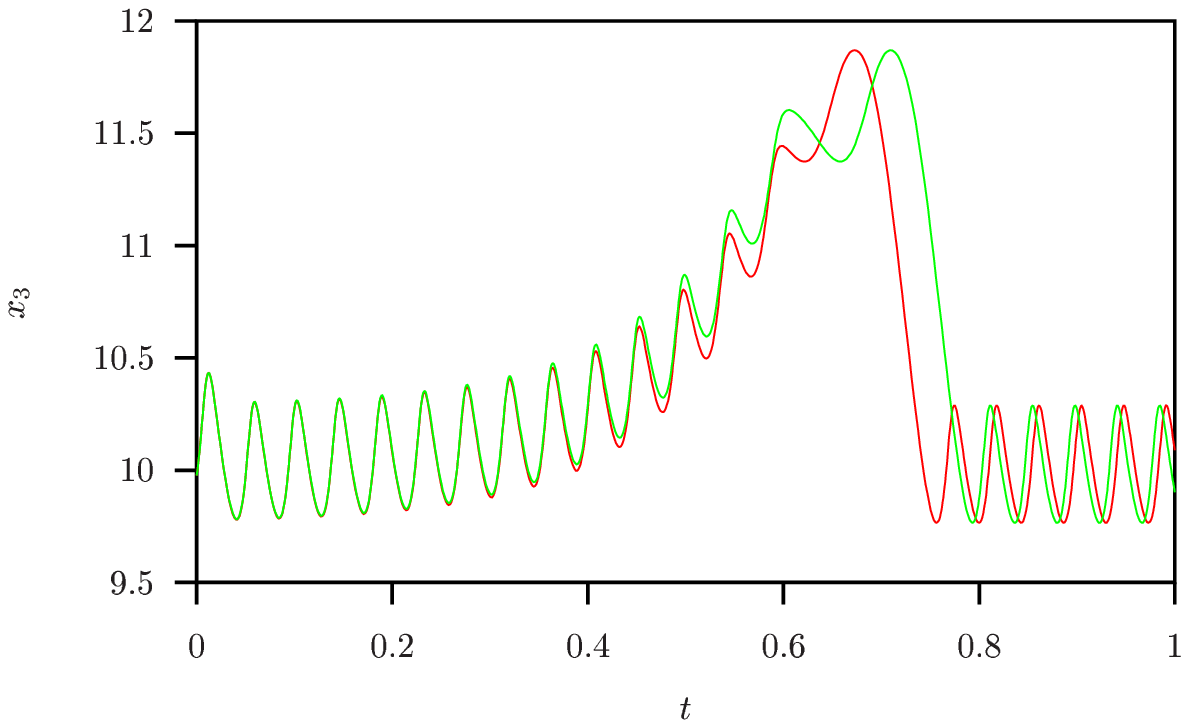}}
    \vspace{0.cm}}\\
  { { \includegraphics{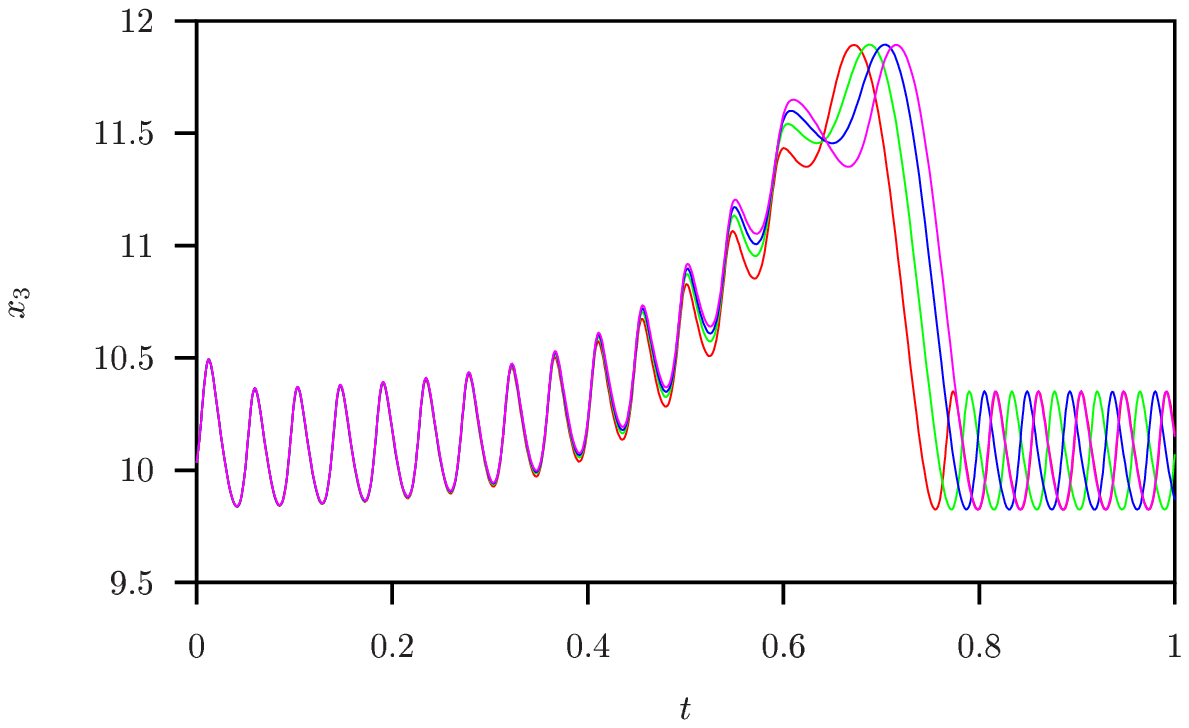}} 
  }}
  \vspace{-0.2cm}
  \caption{
    Profiles of the homoclinic cycle-to-cycle connections for $d_2 =
    0.0125$. The left panel compares the two profiles of the
    connections for $d_1 = 0.2809078$, the right panel compares the
    four profiles for $d_1 = 0.27850$, which is between the limit
    points for the primary and secondary branches. The connection
    time $T$ is scaled to one.}
  \label{fig:cycleprofile}
\end{figure*}

To obtain a good initial connection we consider a BVP like
(\ref{eqn:BVP-homotopy}), with 6 free parameters: $\mu^{\pm}$,
$T^{\pm}$, $T$, and, in turn, one of the four homotopy parameters
$h_{11},h_{12},h_{21},h_{22}$. The selected boundary conditions
(\ref{BVP22c1}) are 
$$
\Phi^{+}[x] = x_2^-(0) - 0.125274 \;,
$$
and 
$$
\Phi^{-}[x] = x_1^+(0) - 0.839789 \;,
$$
so, the first condition uses the $x_2$-coordinate of the initial
base point selected on the cycle, while the second condition uses
the $x_1$-coordinate of the initial base point. Observe, that this selection
is somewhat arbitrary and that one can select other base point coordinates.

In the continuation we want $h_{11} = h_{12} = h_{21} = h_{22} = 0$.
However, there are several solutions, that correspond to connecting
orbits with different numbers of excursions near the limit cycle, both
at the starting and the end-part of the orbit. Observe that the
success of the future continuation in $(d_1,d_2)$ seems to depend highly
on the number of excursions near the cycle at the end-point of the
connecting orbit. In the food chain model a decrease in $d_2$ is
accompanied by a decrease in the numbers of excursions near the cycle
at the end-point of the connection, like a wire around a reel. If this
number is too low, a one-parameter continuation in $d_{1,2}$ will yield
incorrect limit points. Also, two-parameter continuations in $(d_1,d_2)$
will most likely terminate at some point.  Hence a starting orbit is
selected with a sufficient number of excursions near the cycle at the
end-point, with $T = 454.04705$ and $\varepsilon^2 = 0.069414$.

\subsection{Continuation}

The continuation of the connecting orbit can be done in $d_{1,2}$
using the primary BVP~(\ref{eqn:BVP}). Equation~(\ref{BVP-0})
ensures that the base points $x^{\pm}(0) \in O^{\pm}$ become different
(and so do the periodic solutions $x^{\pm}(t)$). 

First, however, using this BVP, the connection can be improved by
increasing the connection time, for the same reason as mentioned above
with regard to the number of excursions near the cycle at the
end-point. The increase in $T$ results in an increase of the number of
excursions near the cycle at the end-point of the connecting orbit.
Then, the continuation in $d_{1,2},T$ results in the detection of four
limit points, of which two are identical.  Figure~\ref{fig:fig8ch4}
shows a continuation in $d_1$ for fixed $d_2 = 0.0125$, that detects
four limit points in one ``run''. Observe that in this way not only
the primary ($d_1 = 0.2809078$, twice, and $d_1 = 0.2305987$), but
also the secondary ($d_1 = 0.2776909$) branch is detected.
Figure~\ref{fig:cycleprofile} shows the profiles of the connecting
orbits for $d_2 = 0.0125$.  Observe that for the region of $0.2776909
< d_1 < 0.2809078$ there are four different connecting orbits with the
same connection time $T$ (see right panel).

\begin{figure*}\hspace{0.cm}
  \scalebox{0.7}{
  { { \includegraphics{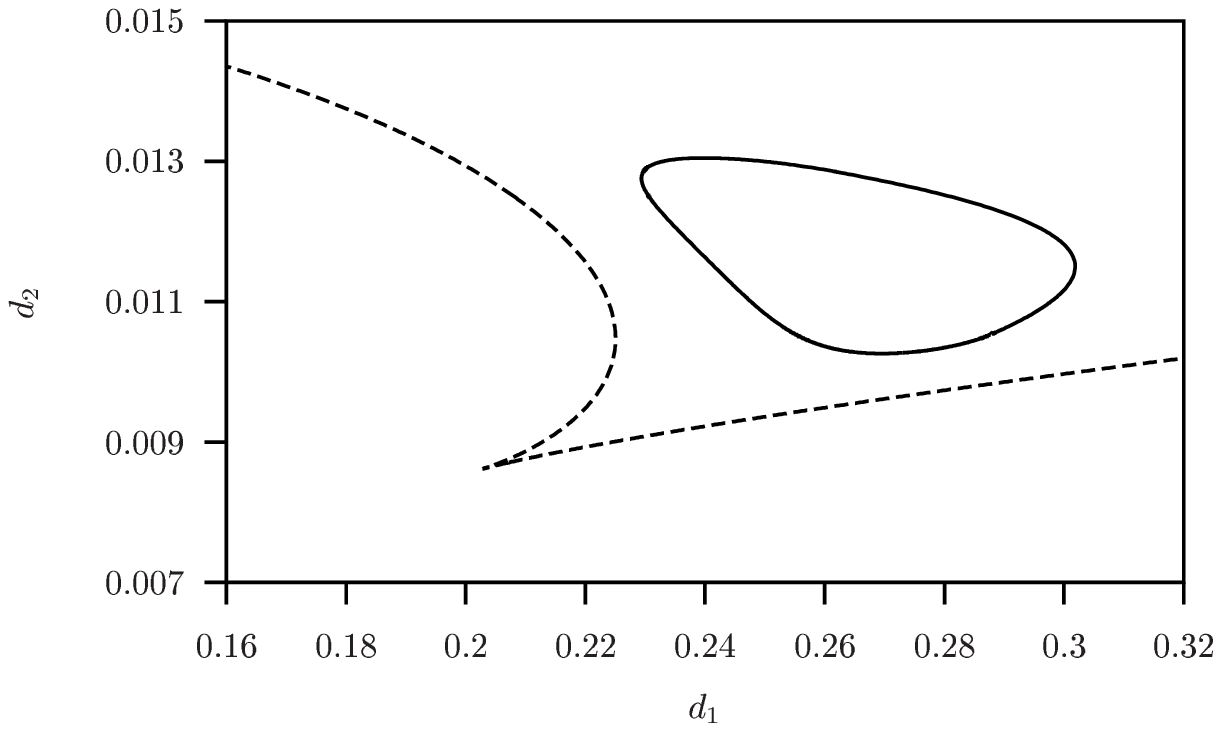}}
    \vspace{0.cm}}\\
  { { \includegraphics{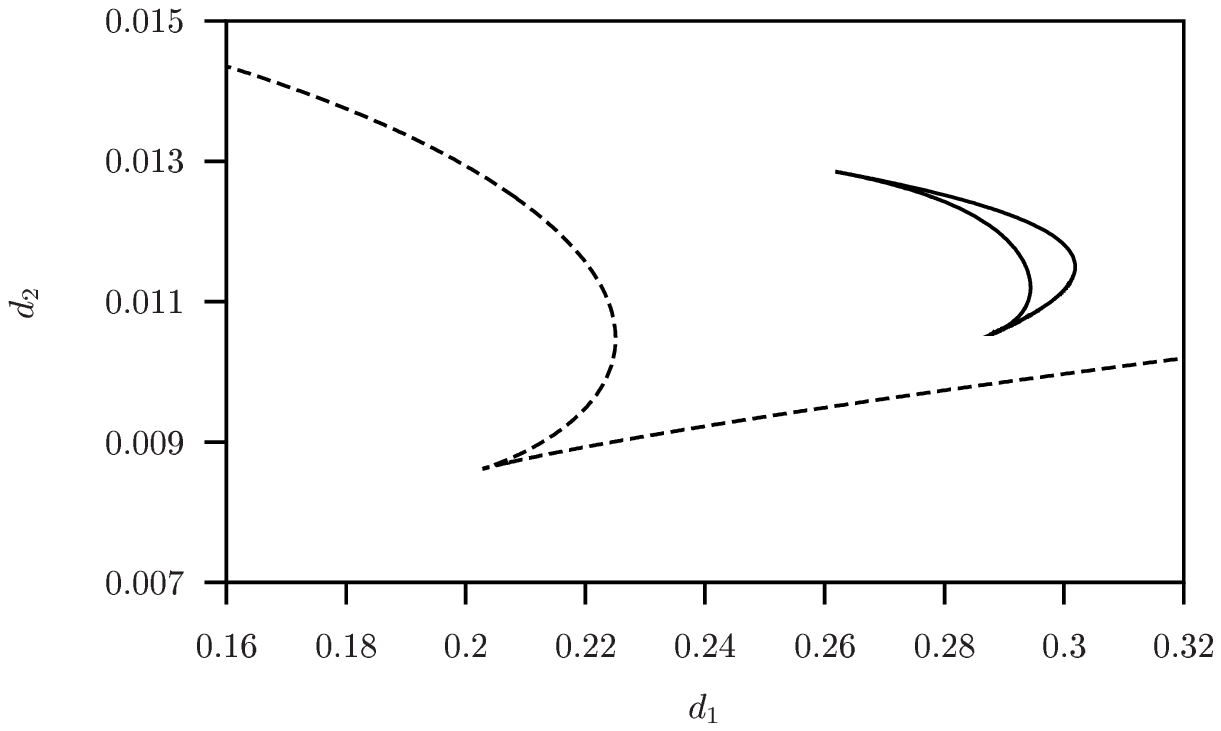}} 
  }}
  \vspace{-0.2cm}
  \caption{
    Two-parameter curves of the primary (left) and secondary (right)
    homoclinic cycle-to-cycle connections of the food chain model.
    Depicted is also the fold bifurcation curve of a limit cycle
    (dashed).}
  \label{fig:connectrm}
\end{figure*}

Using the standard fold-following facilities for BVPs in {\sc auto},
both critical homoclinic orbits can be continued in two parameters
($d_1,d_2$). Along these orbits the stable and unstable invariant
manifold of the cycle are tangent. Starting from $d_1 = 0.2809078$ we
continue the primary branch. The secondary branch is continued from
$d_1 = 0.2776909$. Both curves are depicted in
Figure~\ref{fig:connectrm}.

\section{Discussion}

Our continuation method for cycle-to-cycle connections, using
homotopies in a boundary value setting, is a modified method proposed
in our previous paper for the continuation of point-to-cycle
connections (Doedel et al., 2007\nocite{Doedeletal2007}). The results
discussed here seem to be both robust and time-efficient. Detailed
{\sc auto} demos performing the computations described in Section 6
are freely downloadable from 
\par
\noindent
{\tt www.bio.vu.nl/thb/research/project/globif.}  \par

Provided that the cycle has one simple unstable multiplier, the
proposed method can be extended directly to homoclinic cycle-to-cycle
connections in $n$-dimen\-sio\-nal systems.

\section{Acknowledgements}

The research of the first author (GvV) is supported by the Dutch
Organization for Scientific Research (NWO-CLS) grant no.\ 635,100,013.

\bibliographystyle{ijbc} 
\bibliography{HOM}

\end{document}